\documentclass[10pt]{article}

\usepackage{latexsym,amscd,amsmath,amsthm,amssymb}

\usepackage{mathrsfs} 

\DeclareFontShape{OT1}{cmr}{bx}{sc}
      {
      <-> cmbxsc10
      }{}

\def\proofont{\fontseries{bx}\fontshape{sc}\selectfont}

\renewcommand{\proof}[1][\proofname]{{\par\parindent 0pt {\proofont #1}}}
\newcommand{\QED}{$\hfill\boldsymbol\Box$}

\newcommand{\bbar}{\overline}
\newcommand{\C}{\mathscr{C}}
\newcommand{\HOM}{\mathscr{H}}

\newtheorem{theorem}{Theorem}
\newtheorem{lemma}{Lemma}
\newtheorem{corollary}{Corollary}

\theoremstyle{definition}
\newtheorem{definition}{Definition}

\author{G\'{a}bor Luk\'{a}cs\thanks{I gratefully acknowledge the 
financial support received from York University and Ontario Graduate 
Scholarship Program that enabled me to do this research.}}

\title{On homomorphism spaces of metrizable 
groups\thanks{2000 Mathematics Subject Classification 22A05}}

\begin{document}

\maketitle

\begin{abstract}
For two not necessarily commutative topological groups $G$ and $K$, let
$\HOM(G,K)$ denote the space of all  continuous homomorphisms from $G$ to $K$ 
with the compact-open topology. We  prove that if $G$ is metrizable and $K$ 
is compact  then $\HOM(G,K)$ is a $k$-space. As a consequence we obtain that 
if $D$ is a  dense subgroup of $G$ then $\HOM(D,K)$ is homeomorphic to 
$\HOM(G,K)$, and if $G$ is separable $h$-complete, then the natural map 
$G \rightarrow \C(\HOM(G,K),K)$ is open onto its image.
\end{abstract}


The aim of the present paper is to generalize the result of
Chasco~\cite{Chas2} that for every abelian metrizable group $G$, its dual
group $\hat G$ (i.e. the group of homomorphisms into the unit circle,
$\mathbb{T}$) is a $k$-space under the compact-open topology.  We prove that 
the space of homomorphisms $\HOM(G,K)$ is a $k$-space whenever $G$ is a (not 
necessarily commutative) metrizable topological group and $K$ is a compact 
topological group which satisfies assumptions that we call 
``radical-based" below.

\begin{definition}
A topological group $K$ is {\it radical-based}, if it 
has a countable base $\{\Lambda_n\}$ at $e$,  such that each 
$\Lambda_n$  is symmetric, and for all $n \in \mathbb{N}$: 

\begin{list}{{\rm (\arabic{enumi})}}
{\usecounter{enumi}\setlength{\labelwidth}{25pt}}

\item
$(\Lambda_n)^n \subset \Lambda_1$; 

\item
$a^1, a^2 \ldots, a^n \in  \Lambda_1$  implies $a \in \Lambda_n$.
\end{list}
\end{definition}

Any topological subgroup $K$ of the unitary group of a 
$C^*$-algebra  is radical-based: one can define 
$\Lambda_n = \{ u \in K \mid \| u - e \| <  \varepsilon_n \}$ for a 
suitably chosen sequence $\{\varepsilon_n \}$.

\bigskip

Recall that a Hausdorff topological space $X$ is called a {\it k-space} if 
$F\subset  X$ is closed  if and only if $F \cap C$ is closed for 
every closed compact subset $C$ of $X$.

\begin{theorem} \label{theorem:Hkspace}
For a metrizable topological group $G$ and a radical-based compact group 
$K$, $\HOM(G,K)$ is a $k$-space.
\end{theorem}


In order to prove Theorem~\ref{theorem:Hkspace} we will need the two 
results below. To shorten notations, for $\alpha \in \HOM(G,K)$
we put $S_\alpha(A,B) = S(A,B)\alpha \cap \HOM(G,K)$
where $S(A,B) = \{ \gamma \mid \gamma (A) \subset B\}$ 
($A \subset G$ and $B \subset K$).

\begin{lemma} \label{lemma:lambda}
If $K$ is radical-based then its base $\{\Lambda_n\}$ at $e$ satisfies:

\noindent
{\rm (a)}
$\Lambda_{2k} \Lambda_{2k} \subset  \Lambda_k$ for all $k \in \mathbb{N}$;

\noindent
{\rm (b)}
$\bbar{\Lambda_{2k}} \subset \Lambda_k$ for all $k \in \mathbb{N}$.
\QED
\end{lemma}

\begin{lemma} \label{lemma:Ucompact}
Suppose that $G$ is metrizable and $K$ 
is compact and radical-based. Let $\alpha \in \HOM(G,K)$ and $U$ be a 
neighborhood of  $e$ in $G$. Then $S_\alpha(U,\Lambda_2)$ is precompact, 
i.e. $\bbar{S_\alpha(U,\Lambda_2)}$ is compact. \QED
\end{lemma}

\medskip

\proof[Proof of Theorem~\ref{theorem:Hkspace}.]  Let $\Phi \subset \HOM(G,K)$ 
be a set such that for any compact subset $\Xi$ of $\HOM(G,K)$,
$\Phi \cap \Xi$ is closed. 
We have to  prove that $\Phi$ is closed. To that end let 
$\zeta \in \HOM(G,K)$ such that $\zeta \not \in \Phi$. It suffices to find a 
compact subset $C$ of $G$ and $l \geq 2$ such that 
$S_\zeta (C,\Lambda_{2l}) \cap \Phi = \emptyset$.

The group $G$ is first countable, so let $\{U_n\}_{n=1}^\infty$ be a base 
at $e$. We  may assume that $U_n$ is decreasing. Set $U_0 = G$. We are 
going to find  $l \geq 2$ and construct inductively a family 
$\{F_n\}_{n=0}^\infty$ of  finite subsets of $G$ such that for all  
$n \geq 0$
\begin{eqnarray*}
& \mbox{(1)} &
F_n \subset U_n \mbox{,}  \\
& \mbox{(2)} &
\bigcap\limits_{k=1}^n S_\zeta ( F_k, \bbar{\Lambda_{2l}} ) 
\cap \bbar{ S_\zeta (U_{n+1}, \Lambda_2 )} \cap \Phi = 
\emptyset \mbox{.}
\end{eqnarray*}
First we have to construct $F_0$. By Lemma~\ref{lemma:Ucompact}, 
$\bbar { S_\zeta(U_1,\Lambda_2)}$ is compact, thus by the 
assumption $\bbar { S_\zeta(U_1,\Lambda_2)} \cap \Phi$ is 
closed. On compact subsets of $\C(G,K)$ the compact-open topology 
coincides with the topology of pointwise convergence. But 
$\zeta \not\in \bbar { S_\zeta(U_1,\Lambda_2)} \cap \Phi$, so 
there exists a neighborhood of $\zeta$ in the pointwise topology which is 
disjoint from $\bbar { S_\zeta(U_1,\Lambda_2)} \cap \Phi$.
It is clear that sets of the form $S_\zeta(F,\Lambda_l)$ 
where $F \subset G$ is finite form a base at $\zeta$ for the pointwise 
topology on $\HOM(G,K)$. So there exists $F_0$ such that 
\begin{equation}
S_\zeta(F_0,\Lambda_l) \cap \bbar {S_\zeta(U_1,\Lambda_2)} \cap \Phi = 
\emptyset \mbox{.}
\end{equation}
(Without loss of generality we may assume $l \geq 2$.) 
By Lemma~\ref{lemma:lambda}, $\bbar{\Lambda_{2l}} \subseteq \Lambda_l$, 
thus $S_\zeta(F_0,\bbar{\Lambda_{2l}}) \subset S_\zeta(F_0,\Lambda_l)$. 
In particular:
\begin{equation}
S_\zeta(F_0,\bbar{\Lambda_{2l}}) \cap \bbar{S_\zeta(U_1,\Lambda_2)} \cap \Phi 
= \emptyset \mbox{.}
\end{equation}

Suppose that we have already constructed $F_0,\ldots,F_{n-1}$ such that 
(1) and (2) hold. For all  $x \in U_n$ we define
\begin{equation}
\Delta_x = \bigcap\limits_{k=0}^{n-1} 
S_\zeta(F_k,\bbar{\Lambda_{2l}}) \cap
S_\zeta(\{x\},\bbar{\Lambda_{2l}}) \cap 
\bbar{S_\zeta(U_{n+1},\Lambda_2)} \cap \Phi \mbox{.}
\end{equation}
Notice, that the sets $\Delta_x$ are closed, because each 
$S_\zeta(F_k,\bbar{\Lambda_{2l}})$ is closed even in the 
pointwise  topology. But then
\begin{equation}
\bigcap\limits_{x \in U_n} \Delta_x =
\bigcap\limits_{k=0}^{n-1}  S_\zeta(F_k,\bbar{\Lambda_{2l}}) \cap
S_\zeta(U_n,\bbar{\Lambda_{2l}}) \cap
\bbar{S_\zeta(U_{n+1},\Lambda_2)} \cap \Phi \mbox{.}
\end{equation}
Since $S_\zeta(U_n,\bbar{\Lambda_{2l}}) \subset S_\zeta(U_{n},\Lambda_2)$, 
this means (using assumption (2)) that
\begin{equation}
\bigcap\limits_{x \in U_n} \Delta_x \subset
\bigcap\limits_{k=0}^{n-1}
S_\zeta(F_k,\bbar{\Lambda_{2l}}) \cap
\bbar{S_\zeta(U_{n},\Lambda_2)} \cap \Phi = \emptyset \mbox{.}
\end{equation}
$\Delta_x$ are closed subsets of 
$\bbar{S_\zeta(U_{n+1},\Lambda_2)}$, which is compact by
Lemma~\ref{lemma:Ucompact}. Therefore, there must be a finite set 
$F_n \subset U_n$ such that 
$\bigcap\limits_{x \in F_n} \Delta_x = \emptyset$, in other words:
\begin{equation}
\bigcap\limits_{k=0}^{n-1}  S_\zeta(F_k,\bbar{\Lambda_{2l}}) \cap
S_\zeta(F_n,\bbar{\Lambda_{2l}}) \cap 
\bbar{S_\zeta(U_{n+1},\Lambda_2)} \cap \Phi = \emptyset 
\mbox{,}
\end{equation}
as desired.

Let $C = \bigcup\limits_{n=0}^\infty F_n \cup \{e\}$. We have $F_n \subset 
U_n$, so $C$ is a set of elements converging to $e$. Thus $C$ is sequentially 
compact, but since $G$ is metrizable, it means that $C$ is compact. 
It is clear that 
$S_\zeta(C, \bbar{\Lambda_{2l}}) \cap\bbar{S_\zeta(U_n,\Lambda_2)} \cap \Phi 
= \emptyset$. Since 
$\HOM(G,K) = \bigcup\limits_{n=1}^\infty S_\zeta(U_n,\Lambda_2)$, this means 
that $S_\zeta(C, \bbar{\Lambda_{2l}}) \cap \Phi = \emptyset$. 
Therefore $S_\zeta(C, \Lambda_{2l} ) \cap \Phi = \emptyset$. 
\QED


\bigskip

A topological space $X$ is {\it hemicompact} if 
$X$ is the countable union of compact subspaces $X_n$, such that every  
compact subset of $X$ is contained in a finite union of the sets $X_n$.

\begin{corollary}
For a metrizable topological group $G$ and a radical-based compact group
$K$, $\C(\HOM(G,K),K)$ is completely metrizable.
\end{corollary}

\proof Once metrizability has been shown the completeness is obvious, 
because $\HOM(G,K)$ is a $k$-space, and $K$ is complete (because it is compact). 
Since in \cite{Arens} it was shown that if $X$ is hemicompact then
$\C(X,K)$ is metrizable, it suffices to show that $\HOM(G,K)$ is hemicompact. 

Let $\Xi$ be a compact subset of $\HOM(G,K)$. By the Ascoli Theorem $\Xi$ is 
equicontinuous, in particular  there exists a neighborhood $U$ of $e$ such 
that  $\xi(U) \subset \Lambda_2$ for all $\xi \in \Xi$. In other words, 
$\Xi \subset S_e(U,\Lambda_2)$. Let $\{U_n\}$ be a base at $e \in G$.
For some $n \in \mathbb{N}$, $U_n \subset U$, thus
$\Xi \subset S_e(U,\Lambda_2) \subset S_e(U_n,\Lambda_2)$. By 
Lemma~\ref{lemma:Ucompact}, $\bbar{S_e(U_n,\Lambda_2)}$ is compact, and 
clearly $\HOM(G,K)=\bigcup\limits_{n=1}^\infty \bbar{S_e(U_n,\Lambda_2)}$, 
hence  $\HOM(G,K)$ is hemicompact, as desired. \QED

\begin{definition}
$G$ is {\it h-complete} if for any continuous homomorphism $f: G
\rightarrow H$ the subgroup $f(G)$ is closed in $H$.
\end{definition}

The Corollary below generalizes a theorem by Chasco~\cite{Chas2} stating
that for a separable metrizable complete abelian group $G$, the natural
map $G \rightarrow \hat {\hat G}$ is an isomorphism of topological groups
if and only if it is bijective.

\begin{corollary} \label{corollary:unit}
Suppose that $G$ is separable metrizable and h-complete, and suppose 
further that $K$ is  compact and radical-based. Then, the natural map 
$N: G \longrightarrow \C(\HOM(G,K),K)$ is continuous and 
open onto the image. (In  particular it is an embedding if and only if 
it is one-to-one.)
\end{corollary}

\proof
Since $G$ is metrizable it is clear that $N$ is continuous, because the 
natural map $G \rightarrow \C(\C(G,K),K)$ is continuous.
To see  that it is open onto its image, we notice that 
since $G$ is $h$-complete,  $N(G)$ is closed in $\C(\HOM(G,K),K)$, and 
therefore it is complete metric. Hence $N(G)$ is a Baire space. 
Applying an open map type of theorem (\cite[Cor. 32.4]{Hus}) one obtains
that  $N$ is open onto the image. 
\QED

\begin{theorem} \label{theorem:dense}
Let $K$ be a compact radical-based group.
If $D$ is a dense subgroup of the metrizable group $G$ 
then $\HOM(D,K) \cong \HOM(G,K)$.
\end{theorem}

\proof 
Clearly, $\HOM(D,K) = \HOM(G,K)$ as sets, and we have an induced 
map $\iota: \HOM(G,K) \longrightarrow \HOM(D,K)$ by restriction. Since $\iota$ 
is continuous, and bijective, we only have to show that $\iota$ is open. 
To that end we will show that the inverse image 
of a compact set is compact. Since $\HOM(D,K)$ is a $k$-space it will 
imply that $\iota$ is open.

Take a compact subset $\Phi$ of $\HOM(D,K)$. Then, by the Ascoli Theorem $\Phi$
is equicontinuous; in particular there exists a neighborhood $U$ of $e$ in $G$,
such that $\zeta(U\cap D) \subset \Lambda_4$ for all $\zeta \in \Phi$, and so 
$\zeta(\bbar{U\cap D}) \subset \bbar{\Lambda_4}\subset \Lambda_2$ for all 
$\zeta \in \Phi$. Thus $\Phi \subset S_e (\bbar{U\cap D},\Lambda_2)$.

Let $V$ be a symmetric neighborhood of $e$ in $G$ such that 
$V^2 \subset U$, and let $x \in V$. There exists a sequence 
$\{x_n\} \subset D$  such that $x_n \rightarrow x$, and thus for 
$n \geq n_0$, $x_n \in xV \subset VV \subset U$.
Since $x_n \in D$,  $x_n \in U\cap D$ for 
$n \geq n_0$,  thus $x \in \bbar{U \cap D}$, and hence 
$V \subset \bbar{U \cap D}$.  Therefore
$\Phi \subset S_e (\bbar{U\cap D},\Lambda_2) \subset S_e (V,\Lambda_2)$. 
By Lemma~\ref{lemma:Ucompact}, $S_e (V,\Lambda_2)$ is precompact, hence 
$\Phi$ is compact in $\HOM(G,K)$ as it is closed there. \QED

\bigskip

We note that Theorem~\ref{theorem:dense} is a generalization to the
non-abelian case of a similar result by Chasco~\cite{Chas2}.

\section*{Acknowledgements}

\mbox{ }

The author is deeply indebted to his PhD thesis supervisor,
Prof.~Walter~Tholen, for his dedicated mentorship that made this research
possible. 

\bigskip

He thanks Dr.~Gavin~Seal, for his assistance in simplifying some
proofs and preparing the paper for submission, and he also appreciates the
helpful comments of the anonymous referee. 

\bigskip

Last, but not least, the author offers special thanks to his students for
their enormous encouragement.

{\bigskip\bigskip\noindent Department of Mathematics \& Statistics\\
York University, 4700 Keele Street\\
Toronto, Ontario, M3J 1P3\\
Canada

\bigskip\noindent{\em e-mail: lukacs@mathstat.yorku.ca} }

\end{document}